\documentclass[10pt]{amsart}
\usepackage{latexsym}
\usepackage{amsmath}
\usepackage{amssymb}
\usepackage[all]{xy}

\numberwithin{equation}{section}
\newtheorem{thm}{Theorem }[section]

\newtheorem{lem}[thm]{Lemma}

\newtheorem{rem}[thm]{Remark}

\newcommand{\Z}{{\mathbb Z}}
\newcommand{\R}{{\mathbb R}}
\newcommand{\C}{{\mathbb C}}
\newcommand{\Q}{{\mathbb Q}}

\begin{document}
\title{Examples of rational toral rank complex }
\author{Toshihiro YAMAGUCHI }
\dedicatory{ Dedicated to  Yves F\'{e}lix on his 60th birthday}
\address{Faculty of Education, Kochi University, 2-5-1,Kochi,780-8520, JAPAN}
\email{tyamag@kochi-u.ac.jp}
\maketitle

\begin{abstract}
In \cite[Appendix]{Y}, we see  a CW complex ${\mathcal T}(X)$,
which gives a rational  homotopical  classification of  almost free toral actions on 
spaces in the rational homotopy type of $X$ 
associated  with
rational toral ranks
and also presents certain  relations in them.
We  call it the {\it  rational toral rank complex} of $X$.
It represents a variety of toral actions. 
In this note, we will give  effective   2-dimensional examples of it
when $X$ is a finite product of odd spheres.
This is a combinatorial approach in  rational homotopy theory.   
\end{abstract}

\section{Introduction}
Let $X$ be a simply connected CW complex with $\dim H^*(X;\Q )<\infty$
and  $ r_0(X)$ be the {\it rational  toral rank}
 of   $X$,
 which is the largest integer $r$ such that an $r$-torus
 $T^r=S^1 \times\dots\times S^1$($r$-factors)  can act continuously
 on a CW-complex $Y$ in  the rational homotopy type of $X$
 with all its isotropy subgroups finite (such an action  is called {\it almost free}) \cite{H}.
It is a very  interesting rational invariant.
For example,  
the inequality $$r_0(X)=r_0(X)+r_0(S^{2n})<r_0(X\times S^{2n})\ \ \ \ (*)$$ can hold  for a
formal  space $X$ and an integer $n>1$ \cite{JL}. 
It must appear
 as
one phenomenon
in a variety of almost free toral actions.
 The example $(*)$ is given due to  Halperin 
 by using {\it Sullivan minimal model} \cite{FHT}.
 
Put the Sullivan minimal model $M(X)=(\Lambda V,d)$ of $X$.
If an $r$-torus $T^r$ acts on $X$
by $\mu :T^r\times X\to X$, there is a minimal KS extension
 with $|{t_i}|=2$ for $i=1,\dots,r$
$$
(\Q[t_1,\dots,t_r],0)
 \to (\Q[t_1,\dots,t_r] \otimes \wedge{V},D)
 \to (\wedge{V},d)$$
with  $Dt_i=0$ and
$Dv \equiv dv$ modulo the ideal $(t_1,\dots,t_r)$ for $v\in V$
which is induced from the Borel fibration \cite{FOT}$$
X \to ET^r \times_{T^r}^{\mu} X \to BT^r.
$$
\noindent
According  to \cite[Proposition 4.2]{H},
 $r_0(X) \ge r$ if and only if there is a KS extension of above 
 satisfying $\dim H^*(\Q[t_1,\dots,t_r] \otimes \wedge{V},D)<\infty$.
 Moreover,
  then $T^r$ acts freely on a finite complex 
 that has  the same rational homotopy type as $X$.
So we will discuss this note by Sullivan models.

We want to   give  a classification of rationally almost free toral actions on $X$ associated   with
rational toral ranks
and also present certain  relations in them. 
Recall  a finite based  CW complex ${\mathcal T}(X)$ in \cite[\S 5]{Y}.
Put ${\mathcal X}_r=\{ (\Q[t_1,\dots,t_r] \otimes \wedge{V},D) \}$
the set of isomorphism classes  of KS-extensions of $M(X)=(\Lambda V,d)$
such that  $\dim H^*(\Q[t_1,\dots,t_r] \otimes \wedge{V},D) <\infty$.
First, the set of 0-cells
${\mathcal T}_0(X)$ is the finite sets $\{ (s,r)\in \Z_{\geq 0} \times \Z_{\geq 0} \}$
where the point $P_{s,r}$ of the coordinate $(s,r)$ exists if there is a model 
$(\Lambda W,d_W)\in {\mathcal X}_r$ 
and $r_0(\Lambda W,d_W)=r_0(X)-s-r$.
Of course, the model may  not be uniquely determined.
Note the base point $P_{0,0}=(0,0)$ always exists by $X$ itself.

Next, 1-skeltons (vertexes)   of   the   1-skelton  ${\mathcal T}_1(X)$
are  represented  by a KS-extension $(\Q [t],0)\to (\Q [t]\otimes \Lambda W,D)\to (\Lambda W,d_W)$
with  $\dim H^*(\Q[t] \otimes \wedge{W},D) <\infty$
for  $(\Lambda W,d_W) \in{\mathcal X}_r$, 
where $W=\Q( t_1,\dots,t_r) \oplus {V}$
and $d_W|_V=d$.
It is given as 
{\small 
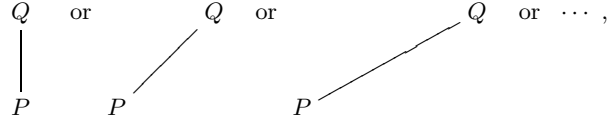
$$\xymatrix{ 
Q\\
P \ar@{-}[u]
}\ \ \ \ \mbox{or} \ 
\xymatrix{ 
&Q\\
P \ar@{-}[ru]&
}
\ \ \  \mbox{or} \ \xymatrix{ 
&&Q\\
P \ar@{-}[rru]&&
} \ \mbox{\ \ or} \ \ \cdots \ ,
$$
}
\noindent
where $P$ exists by  $(\Lambda W,d_W)$ and $Q$ exists by  $(\Q [t]\otimes \Lambda W,D)$.
The 2-cell is given if there is a (homotopy) commutative diagram of restrictions 
 {\small 
$$\xymatrix{ (\Lambda W,d_W)\ \  &
(\Q [t_{r+2}]\otimes \Lambda W,D_{r+2})\ar[l]\\ 
(\Q [t_{r+1}]\otimes \Lambda W,D_{r+1})\ar[u]&(\Q [t_{r+1},t_{r+2}]\otimes \Lambda W,D), \ar[u]\ar[l]& }$$
}
\noindent
which represents
(a horizontal deformation of) 
{\small
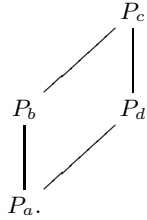
$$\xymatrix{ 
&  P_c \\
P_b \ar@{-}[ur] & P_d\ar@{-}[u]\\
P_a.\ar@{-}[u]\ar@{-}[ur]&\\
}$$
}
\noindent
Here $P_a$  exists by $(\Lambda W,d_W)$,
$P_b$(or $P_d$)  by $(\Q [t_{r+1}]\otimes \Lambda W,D_{r+1})$,
$P_c$ by $(\Q [t_{r+1},t_{r+2}]\otimes \Lambda W,D)$
and 
$P_d$(or $P_b$)   by $(\Q [t_{r+2}]\otimes \Lambda W,D_{r+2})$.
Then we say that a 2-cell attachs to (the tetragon)  $P_aP_bP_cP_d$.
Thus we  can  construct the 2-skelton  ${\mathcal T}_2(X)$.



Generally, an  n-cell is 
given by an $n$-cube where
a vertex of $(\Q [t_{r+1},..,t_{r+n}]\otimes \Lambda W,D)$
 of height $r+n$,
 n-vertexes
$\{ (\Q [t_{r+1},..,\overset{\vee}{t_{r+i}},..., t_{r+n}]\otimes \Lambda W,D_{(i)})\}_{1\leq i\leq n}$
 of height $r+n-1$, 
 $\cdots$ ,
 a vertex $( \Lambda W, d_W)$
  of height $r$.
  Here $\vee$ is the symbol which removes
the below element and  the differential $D_{(i)}$
is the restriction of $D$.

We will call this connected  regular
complex  ${\mathcal T}(X)=\cup_{n\geq 0}{\mathcal T}_n(X)$
the {\it rational toral rank complex} (r.t.r.c.) of $X$.
Since $r_0(X)<\infty$ in our case, it is a finite complex.
For example, 
when $X=S^3\times  S^3$ and $Y=S^5$,
we have 
$${\mathcal T}(X)\vee {\mathcal T}(Y)=
{\mathcal T}_1(X)\vee {\mathcal T}_1(Y)={\mathcal T}_1(X\times Y)=
{\mathcal T}(X\times Y),$$
which is an unusual
case. Then, of course, $r_0(X)+r_0(Y)=r_0(X\times Y)$.
Recall that 
 $r_0(S^3\times S^3)+r_0(S^7)=r_0(S^3\times S^3\times S^7)$
but ${\mathcal T}_1(S^3\times S^3)\vee {\mathcal T}_1(S^7)\subsetneq {\mathcal T}_1(S^3\times S^3\times S^7)$ \cite[Example 3.5]{Y}.
In \S 2, we see that r.t.r.c.  is 
not complicated as a CW complex but delicate.
We see in Theorems 2.2 and  2.3 that the differences
between $X= Z\times S^7$ and  $Y= Z\times S^9$ for some products $Z$
of odd spheres
make certain different homotopy types of 
r.t.r.c.,
respectively.
Remark that 
the above inequality $(*)$ is, if anything,  a property on ${\mathcal T}_0(X)$ or  ${\mathcal T}_1(X)$
as the  example of  Theorem \ref{3}(1).
We see in Theorem \ref{3}(2)
an  example that  ${\mathcal T}_1(X)={\mathcal T}_1(X\times \C P^n)$
but
 $ {\mathcal T}_2(X)\subsetneq  {\mathcal T}_2(X\times \C P^n)$,
which is,  in a sense, a  higher dimensional phenomenon of $(*)$.

%





\section{Examples }

In this section, the symbol
$P_iP_jP_kP_l$
means the tetragon,
which is the cycle  with vertexes $P_i$, $P_j$, $P_k$,  $P_l$
and edges $P_iP_j$, $P_jP_k$, $P_kP_l$,  $P_lP_i$.

In general, it is difficult to show that a point of  ${\mathcal T}_0(X)$
does not exist on a certain coordinate.
So the following lemma is useful for our purpose.

\begin{lem}
If $X$ has the rational homotopy type of the product of  finite odd spheres
and finite complex projective spaces,
then $(1,r)\not\in {\mathcal T}_0(X)$ for any $r$.
\end{lem}
 
 \noindent{\it Proof.} 
Suppose that  $X$  has the rational homotopy type of the product of  n odd spheres
and m complex projective spaces.
Put  a minimal model  
$A=(\Q [t_1,\cdots ,t_{n-1},x_1,\cdots ,x_m]\otimes \Lambda (v_1,\cdots ,v_n,y_1,
\cdots ,y_m ),D)$
with $|t_1|=\cdots =|t_{n-1}|=|x_1|=\cdots =|x_m|=2$ and $|v_i|, |y_i|$ odd.
If $\dim H^*(A)<\infty$,
then $A$  
is pure; i.e., $Dv_i, Dy_i\in \Q [t_1,..,t_{n-1},x_1,\cdots ,x_m]$ for all $i$.
Therefore,  from \cite[Lemma 2.12]{JL}, $r_0(A)=1$.
Thus we have $(1,r_0(X)-1)=(1,n-1)\not\in  {\mathcal T}_0(X)$.
 \hfill\qed\\

\begin{thm}
 Put $X=S^3\times S^3\times S^3\times S^{7}\times S^{7}$
and $Y=S^3\times S^3\times S^3\times S^{7}\times S^{9}$.
Then
${\mathcal T}_1(X)={\mathcal T}_1(Y)$.
But  ${\mathcal T}(X)$ is contractible and ${\mathcal T}(Y)\simeq S^2$.
\end{thm}

\noindent
{\it Proof.}
 Let $M(X)=(\Lambda V,0)=(\Lambda (v_1,v_2,v_3,v_4,v_5),0)$ with $|v_1|=|v_2|=|v_3|=3$ and  $|v_4|=|v_5|=7$.
 Then $${\mathcal T}_0(X)=\{ P_{0,0},  P_{0,1}, P_{0,2}, P_{0,3}, P_{0,4}, P_{0,5}, P_{2,1}, P_{2,2}, P_{2,3}, P_{3,1}, P_{3,2} \}.$$
For example, they are given as follows.

(0) $P_{0,0}$ is given by $( \Lambda V,0)$.

(1) $P_{0,1}$ is given by $(\Q [t_1]\otimes \Lambda V,D)$ with 
$Dv_1=t_1^2$ and $Dv_2=Dv_3=Dv_4=Dv_5=0$.  

(2) $P_{0,2}$ is given by $(\Q [t_1,t_2]\otimes \Lambda V,D)$  with 
$Dv_1=t_1^2$, $Dv_2=t_2^2$ and $Dv_3=Dv_4=Dv_5=0$.  

(3) $P_{0,3}$ is given by $(\Q [t_1,t_2,t_3]\otimes \Lambda V,D)$  with 
$Dv_1=t_1^2$, $Dv_2=t_2^2$, $Dv_3=t_3^2$  and $Dv_4=Dv_5=0$.  

(4) $P_{0,4}$ is given by $(\Q [t_1,t_2,t_3,t_4]\otimes \Lambda V,D)$  with 
$Dv_1=t_1^2$, $Dv_2=t_2^2$, $Dv_3=t_3^2$, $Dv_4=t_4^4$  and $Dv_5=0$.  

(5) $P_{0,5}$ is given by $(\Q [t_1,t_2,t_3,t_4,t_5]\otimes \Lambda V,D)$  with 
$Dv_1=t_1^2$, $Dv_2=t_2^2$, $Dv_3=t_3^2$, $Dv_4=t_4^4$  and $Dv_5=t_5^4$.
  
(6) $P_{2,1}$ is given by $(\Q [t_1]\otimes \Lambda V,D)$ with  $Dv_1=Dv_2=Dv_3=Dv_5=0$ and 
$Dv_4=v_1v_2t_1+t_1^4$  

(7) $P_{2,2}$ is given by $(\Q [t_1,t_2]\otimes \Lambda V,D)$  with 
$Dv_1=Dv_2=0$, $Dv_3=t_2^2$, $Dv_4=v_1v_2t_1+t_1^2$ and $Dv_5=0$.  

(8) $P_{2,3}$ is given by $(\Q [t_1,t_2,t_3]\otimes \Lambda V,D)$  with 
$Dv_1=Dv_2=0$,  $Dv_3=t_2^2$, $Dv_4=t_1^2+v_1v_2t_1$ and $Dv_5=t_3^4$.  

(9) $P_{3,1}$ is given by $(\Q [t_1]\otimes \Lambda V,D)$  with $Dv_1=Dv_2=Dv_3=0$, $Dv_4=v_1v_2t_1+t_1^4$ and $Dv_5=v_1v_3t_1$. 

(10) $P_{3,2}$ is given by $(\Q [t_1,t_2]\otimes \Lambda V,D)$  with 
$Dv_4=v_1v_2t_1+t_1^4$ and $Dv_5=v_1v_3t_1+t_2^4$.

(11)  $P_{4,1}$, i.e., a point of the coordinate $(4,1)$
 does  not exist.
Indeed, if it exists, it must be given by  a model 
$(\Q [t_1]\otimes \Lambda V,D)$ 
whose differential is  $Dv_1=Dv_2=Dv_3=0$ and 
$Dv_4,Dv_5\in \Q [t_1]\otimes \Lambda (v_1,v_2,v_3)$
by degree reason.
But, for any $D$
satisfying such conditions,   we have $\dim H^*(\Q [t_1,t_2]\otimes \Lambda V,\tilde{D})<\infty$ for a KS extension 
$$(\Q [t_2],0)\to (\Q [t_1,t_2]\otimes \Lambda V,\tilde{D})\to (\Q [t_1]\otimes \Lambda V,{D}),$$
that is,
$r_0(\Q [t_1]\otimes \Lambda V,D)>0$. 
It contradicts
 the definition of  $P_{4,1}$.

  ${\mathcal T}_1(X)$ is given as
{\small 
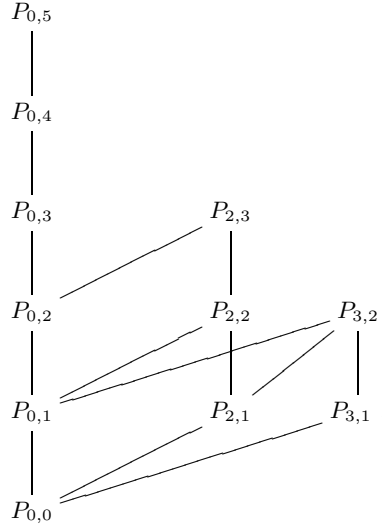
$${\small 
 \xymatrix{ 
P_{0,5}&&&&\\
P_{0,4} \ar@{-}[u]& &  & &\\
P_{0,3} \ar@{-}[u]&  &P_{2,3}&  & \\
P_{0,2} \ar@{-}[u]\ar@{-}[urr]&  & P_{2,2}\ar@{-}[u]&P_{3,2}& \\
P_{0,1} \ar@{-}[u]\ar@{-}[urr] \ar@{-}[urrr]& & P_{2,1}\ar@{-}[u]\ar@{-}[ur]& P_{3,1}\  .\ar@{-}[u] & \\
P_{0,0}\ar@{-}[u]\ar@{-}[urrr]\ar@{-}[urr]\\
}}$$
}

For example, the edges(1-simplexes)  
$$\{ \  P_{0,0}P_{0,1},\ P_{0,1}P_{0,2},\ P_{0,2}P_{0,3},\ P_{0,3}P_{0,4},\ \cdots ,\ P_{0,0}P_{3,1},\  P_{3,1}P_{3,2}\ \} $$
are given as follows.

(1) $P_{0,1}P_{3,2}$ is given by the projection $(\Q [t_1,t_2]\otimes \Lambda V,D)
\to (\Q [t_1]\otimes \Lambda V,{D}_1)$
where  $Dv_1=Dv_2=Dv_3=0$, 
 $Dv_4=v_1v_2t_2+t_1^4$, $Dv_5=v_1v_3t_2+t_2^4$
and  $D_1v_1=D_1v_2=D_1v_3=D_1v_5=0$
and
 $D_1v_4=t_1^4$.

(2) $P_{2,1}P_{3,2}$  is given by 
  $Dv_1=Dv_2=Dv_3=0$, $Dv_4=v_1v_2t_1+t_1^4$ and $Dv_5=v_1v_3t_2+t_2^4$.

(3) $P_{3,1}P_{3,2}$  is given by  $Dv_1=Dv_2=Dv_3=0$, $Dv_4=v_1v_2t_1+t_1^4$
and  $Dv_5=v_1v_3t_1+t_2^4$.

${\bf \cdots}$

 ${\mathcal T}_2(X)$ is given as follows.

(1) $P_{0,0}P_{2,1}P_{3,2}P_{3,1}$ is attached by   a 2-cell from
 $Dv_1=Dv_2=Dv_3=0$,  $Dv_4=v_1v_2(t_1+t_2)+t_1^4$
and  $Dv_5=v_1v_3t_2+t_2^4$.
 (Then $P_{2,1}$ is given by $D_1v_4=v_1v_2t_1+t_1^4$, $D_1v_5=0$
and $P_{3,1}$ is given by $D_2v_4=v_1v_2t_2$, $D_2v_5=v_1v_3t_2+t_2^4$.)
 
(2)  $P_{0,0}P_{0,1}P_{3,2}P_{3,1}$  is attached by   a 2-cell  from
 $Dv_1=Dv_2=Dv_3=0$,  $Dv_4=v_1v_2t_2+t_1^4$
and  $Dv_5=v_1v_3t_2+t_2^4$. 

(3) $P_{0,0}P_{0,1}P_{2,2}P_{2,1}$  is attached by   a 2-cell  from
 $Dv_1=Dv_2=Dv_3=0$,  $Dv_4=v_1v_2t_2+t_2^4$
and  $Dv_5=t_1^4$.

(4) $P_{0,1}P_{0,2}P_{2,3}P_{2,2}$  is attached by   a 2-cell  from
 $Dv_1=Dv_2=0$,  $Dv_3=t_3^2$, $Dv_4=v_1v_2t_2+t_2^4$
and  $Dv_5=t_1^4$.

(5) $P_{0,0}P_{0,1}P_{3,2}P_{2,1}$  is {\it not} attached by  a 2-cell.
Indeed, assume that   a 2-cell attachs on it.
Notice that 
$P_{3,2}$ is given by $(\Q [t_1,t_2]\otimes \Lambda V,D)$ with $Dv_1=Dv_2=Dv_3=0$ and 
$$Dv_4=\alpha (v_1,v_2,v_3)+f,\ \ \ Dv_5=\beta (v_1,v_2,v_3)+g$$
where $\alpha ,\beta \in (v_1,v_2,v_3)$ and 
$\{  f,g\}$ is a regular sequence in $\Q [t_1,t_2]$.
Since  $P_{0,1}P_{3,2}\in {\mathcal T}_1(X)$,
both $\alpha$ and $\beta$ must be  contained in the ideal $(t_i)$ for some $i$.
Also they are  not in $(t_1t_2)$  by degree reason.
Furthermore, since  $P_{2,1}P_{3,2}\in {\mathcal T}_1(X)$, we can put  
that
both $\alpha$ and $\beta$ are contained in 
the monogenetic ideal $(v_iv_j)$ for some $1\leq i<j\leq 3$
without losing generality.
Then,
$\dim H^*(\Q [t_1,t_2,t_3]\otimes \Lambda V,\tilde{D})<\infty$
for a KS extension 
$$(\Q [t_3],0)\to (\Q [t_1,t_2,t_3]\otimes \Lambda V,\tilde{D})\to (\Q [t_1,t_2]\otimes \Lambda V,{D}),$$
by putting $\tilde{D}v_k=t_3^2$ for 
$k\in \{ 1,2,3\}$ with $k\neq i,j$ and $\tilde{D}v_n=Dv_n$ for  $n\neq k$.
Thus we have $r_0(\Q [t_1,t_2]\otimes \Lambda V,D)>0$.
It contradicts to the definition of $P_{3,2}$.

Notice there is no 3-cell
 since it  must attach to a 3-cube (in graphs) in general.
 Thus we see that ${\mathcal T}(X)={\mathcal T}_2(X)$ is contractible.

On the other hand,  let $M(Y)=(\Lambda W,0)=(\Lambda (w_1,w_2,w_3,w_4,w_5),0)$ with $|w_1|=|w_2|=|w_3|=3$,  $|w_4|=7$ and $|w_5|=9$.
Then we see ${\mathcal T}_1(X)={\mathcal T}_1(Y)$
from same arguments.
But, in ${\mathcal T}_2(Y)$, 
$P_{0,0}P_{0,1}P_{3,2}P_{2,1}$ is attached by a 2-cell since we can put 
$Dw_1=Dw_2=Dw_3=0$ and 
$$Dw_4=w_1w_2t_2+t_2^4,\ \ \ Dw_5=w_1w_3t_1t_2+t_1^5,$$
by degree reason.
Here $P_{0,1}$ is given by $D_1w_4=0$, $D_1w_5=t_1^5$
and $P_{2,1}$ is given by $D_2w_4=w_1w_2t_2+t_2^4$,@$D_2w_5=0$.
Others are same as  ${\mathcal T}_2(X)$.
Then three 2-cells on $P_{0,0}P_{0,1}P_{3,2}P_{2,1}$, $P_{0,0}P_{2,1}P_{3,2}P_{3,1}$
and  $P_{0,0}P_{0,1}P_{3,2}P_{3,1}$
in ${\mathcal T}_2(Y)$
makes
{\small
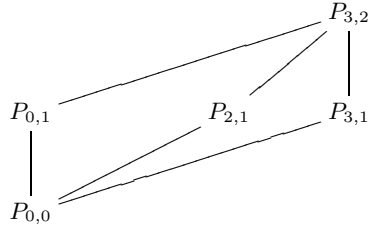
$${\small
 \xymatrix{ 
&  &&P_{3,2}& \\
P_{0,1} \ar@{-}[urrr]& & P_{2,1}\ar@{-}[ur]& P_{3,1}\ar@{-}[u] & \\
P_{0,0}\ar@{-}[u]\ar@{-}[urrr]\ar@{-}[urr]\\
}}$$
}
to be  homeomorphic to $S^2$.
Thus ${\mathcal T}(Y)={\mathcal T}_2(Y)\simeq S^2$.
\hfill\qed\\

\begin{thm}Put $X=S^3\times S^3\times S^3\times S^3\times S^{7}\times S^{7}$
and $Y=S^3\times S^3\times S^3\times S^3\times S^{7}\times S^{9}$.
Then
${\mathcal T}_1(X)={\mathcal T}_1(Y)$.
But  ${\mathcal T}(X)\simeq S^2$ and ${\mathcal T}(Y)\simeq \vee_{i=1}^6S^2_i$.
 \end{thm}

\noindent
{\it Proof}. We see as the proof of Theorem 2.2 that
  $${\mathcal T}_0(X)=\{ P_{0,0},  P_{0,1}, P_{0,2}, P_{0,3}, P_{0,4}, P_{0,5}, P_{0,6}, 
P_{2,1}, P_{2,2}, P_{2,3}, P_{2,4}, $$
$$P_{3,1}, P_{3,2}, P_{3,3}, P_{4,1}, P_{4,2} \}$$
and both 
${\mathcal T}_1(X)$ and ${\mathcal T}_1(Y)$
are  given as 
{\small
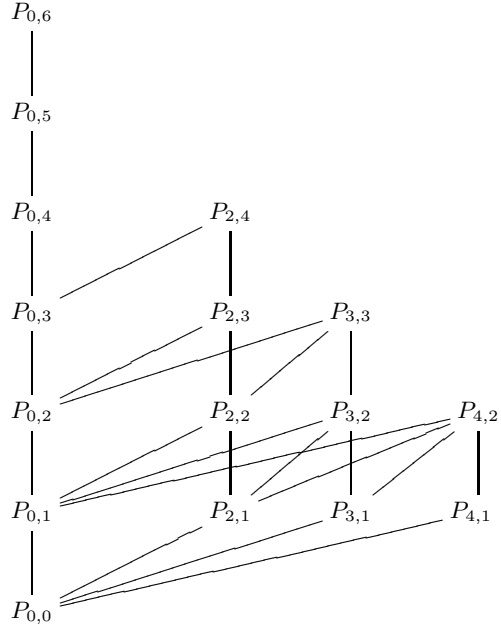
$$
\xymatrix{
P_{0,6}&&& &&\\
P_{0,5} \ar@{-}[u]&&&&&\\
P_{0,4} \ar@{-}[u]&  & P_{2,4} &&& \\
P_{0,3} \ar@{-}[u]\ar@{-}[rru]& 
 & P_{2,3} \ar@{-}[u] & P_{3,3} &&\\
P_{0,2} \ar@{-}[u]\ar@{-}[rru]\ar@{-}[rrru]& & P_{2,2} \ar@{-}[u]\ar@{-}[ru]&
 P_{3,2} \ar@{-}[u]&P_{4,2}&\\
P_{0,1} \ar@{-}[u]\ar@{-}[urr] \ar@{-}[urrr]\ar@{-}[urrrr]&  & P_{2,1}\ar@{-}[u]
\ar@{-}[ru]\ar@{-}[urr]& P_{3,1}  \ar@{-}[u]\ar@{-}[ur]&P_{4,1}\ .\ar@{-}[u]&\\
P_{0,0}\ar@{-}[u]\ar@{-}[urr]\ar@{-}[urrr]\ar@{-}[urrrr]&&\\
}$$
}

 For all tetragons in ${\mathcal T}_1(X)$
except the following 4 tetragons: $$(1)\ P_{0,0}P_{0,1}P_{3,2}P_{2,1}\ 
 \ (2)\ P_{0,1}P_{0,2}P_{3,3}P_{2,2}\ \ (3)\ P_{0,0}P_{0,1}P_{4,2}P_{2,1}\ \ (4)\   P_{0,0}P_{0,1}P_{4,2}P_{3,1},$$
\noindent
2-cells attach  in ${\mathcal T}_2(X)$.
The proof is similar to  it of Theorem 2.2.
Thus we see that 
${\mathcal T}_2(X)$ is homotopy equivalent to 
{\small
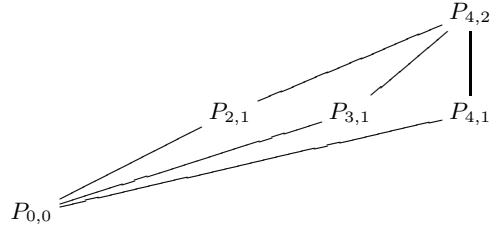
$${\small
 \xymatrix{ 
&  &&&P_{4,2}& \\
&&P_{2,1} \ar@{-}[urr] & P_{3,1}\ar@{-}[ur]& P_{4,1}\ar@{-}[u] & \\
P_{0,0}\ar@{-}[rru]\ar@{-}[urrr]\ar@{-}[urrrr]\\
}}$$
}

\noindent
, which is homeomorphic to $S^2$.
For example, when $M(X)=(\Lambda V,0)=(\Lambda (v_1,v_2,v_3,v_4,v_5,v_6),0)$ with $|v_1|=|v_2|=|v_3|=|v_4|=3$ and  $|v_5|=|v_6|=7$,
2-cells attach 
 $P_{0,0}P_{2,1}P_{4,2}P_{3,1}$, $P_{0,0}P_{3,1}P_{4,2}P_{4,1}$
and $P_{0,0}P_{2,1}P_{4,2}P_{4,1}$ 
from $Dv_1=\cdots =Dv_4=0$,

$Dv_5=v_1v_2t_1+t_1^4,\ 
Dv_6=v_1v_3t_1+v_2v_4t_2+t_2^4$,
 
$Dv_5=v_1v_2t_1+t_1^4, \ Dv_6=v_1v_3(t_1+t_2)+v_2v_4t_2+t_2^4$ \ \ and

$Dv_5=v_1v_2t_1+t_1^4, \ 
Dv_6=v_1v_3t_2+v_2v_4t_2+t_2^4$,\\
respectively.

In ${\mathcal T}_2(Y)$,
2-cells attach all tetragons in 
${\mathcal T}_1(Y)$ by degree reason.
For example,
when
 $M(Y)=(\Lambda W,0)=(\Lambda (w_1,w_2,w_3,w_4,w_5,w_6),0)$ with 
$|w_1|=|w_2|=|w_3|=|w_4|=3$, $|w_5|=7$ and $|w_6|=9$,
put $Dw_1=Dw_2=Dw_3=0$ and 

$(1)\ Dw_4=0,\ Dw_5=w_1w_3t_2+t_2^4,\
Dw_6=w_2w_3t_1t_2+t_1^5$

$(2)\ Dw_4=t_3^2,\ Dw_5=w_1w_3t_2+t_2^4,\
Dw_6=w_2w_3t_1t_2+t_1^5$

$(3)\  Dw_4=0,\ Dw_5=w_1w_2t_2+t_2^4,\ 
Dw_6=w_3w_4t_1t_2+t_1^5$

$(4)\  Dw_4=0,\ Dw_5=w_1w_3t_2+t_2^4,\
Dw_6=w_1w_4t_2^2+w_2w_3t_1t_2+t_1^5$

\noindent
for $(1)\sim (4)$ of above.
Then we can check that ${\mathcal T}(Y)\simeq \vee_{i=1}^6S^2_i$
(${\mathcal T}(Y)$ can not be embedded in $\R^3$).
\hfill\qed\\

\begin{thm}\label{3}
Even when  $r_0(X)=r_0(X\times \C P^n)$ for the $n$-dimensional complex projective space
$\C P^n$,
it does not fold that  ${\mathcal  T}(X)= {\mathcal  T}(X\times \C P^n)$ in general. 
For example,

(1) When $X=S^3\times S^3\times S^3\times S^3\times S^7$ and $n=4$, then 
${\mathcal  T}_1(X)\subsetneq {\mathcal  T}_1(X\times \C P^4)$.

(2) When $X=S^3\times S^3\times S^3\times S^{7}\times S^{7}$ and $n=4$,  then 
${\mathcal  T}_1(X)={\mathcal  T}_1(X\times \C P^4)$
but ${\mathcal  T}_2(X)\subsetneq {\mathcal  T}_2(X\times \C P^4)$.
\end{thm}

\noindent
{\it Proof}.
Put $M( \C P^n)=(\Lambda (x,y),d)$
with $dx=0$ and  $dy=x^{n+1}$ for  $|x|=2$ and $|y|=2n+1$.
Put 
$(\Q [t_1,..,t_r]\otimes \Lambda V\otimes \Lambda (x,y),D)$
the model  of a Borel space 
$ET^r\times_{T^r}(X\times \C P^n)$
of  $X\times \C P^n$.



(1) ${\mathcal  T}_1(X)$ and ${\mathcal  T}_1(X\times \C P^4)$ are given as

{\small
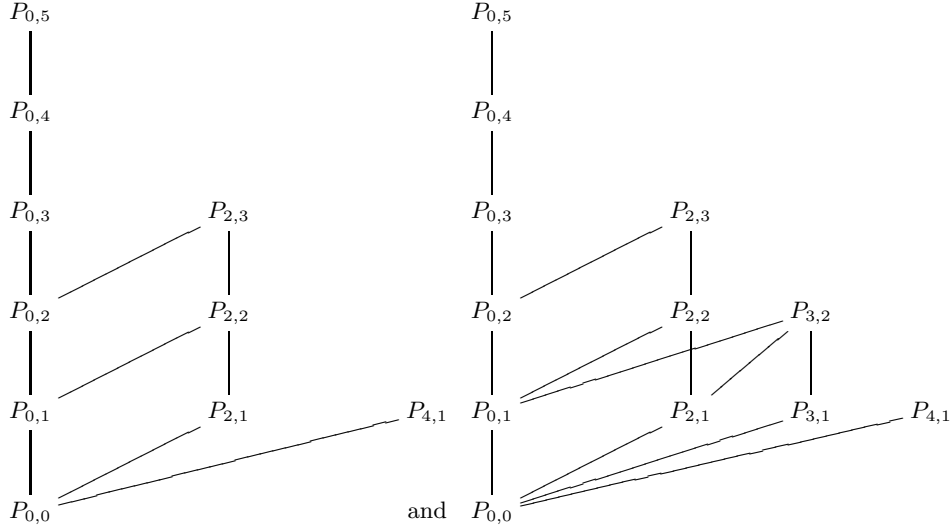
$${\small
 \xymatrix{ 
P_{0,5}&&&&\\
P_{0,4} \ar@{-}[u]& &  & &\\
P_{0,3} \ar@{-}[u]&  &P_{2,3}&  & \\
P_{0,2} \ar@{-}[u]\ar@{-}[urr]&  & P_{2,2}\ar@{-}[u]&& \\
P_{0,1} \ar@{-}[u]\ar@{-}[urr] & & P_{2,1}\ar@{-}[u]& &  P_{4,1}\\
P_{0,0}\ar@{-}[u]\ar@{-}[urr]\ar@{-}[urrrr]&&&&{\rm and}\\
}}\mbox{  }
{\small
 \xymatrix{ 
P_{0,5}&&&&\\
P_{0,4} \ar@{-}[u]& &  & &\\
P_{0,3} \ar@{-}[u]&  &P_{2,3}&  & \\
P_{0,2} \ar@{-}[u]\ar@{-}[urr]&  & P_{2,2}\ar@{-}[u]&P_{3,2}& \\
P_{0,1} \ar@{-}[u]\ar@{-}[urr] \ar@{-}[urrr]& & P_{2,1}\ar@{-}[u]\ar@{-}[ur]& P_{3,1}\ar@{-}[u] & P_{4,1}\\
P_{0,0}\ar@{-}[u]\ar@{-}[urrr]\ar@{-}[urr]\ar@{-}[urrrr]\\
}}$$
}
\noindent
, respectively.
For $M(X)=(\Lambda V,0)=(\Lambda (v_1,v_2,v_3,v_4,v_5),0)$ with $|v_1|=|v_2|=|v_3|=|v_4|=3$ and  $|v_5|=7$.
Here $P_{4,1}$ is given by $Dv_i=0$ for $i=1,2,3,4$ and 
$Dv_5=v_1v_2t_1+v_3v_4t_1+t_1^4$.
It is contained  in both ${\mathcal T}_0(X)$
and ${\mathcal T}_0(X\times \C P^4)$.
On the other hand, $P_{3,2}$ is given by $Dv_i=0$ for $i=1,2,3$,  
$Dv_4=t_2^2$, $Dv_5=v_1v_2t_1+t_1^4$, $Dx=0$ and $Dy=x^5+v_1v_3t_1^2$.
Then  $P_{3,1}$ is given by $Dv_i=0$ for $i=1,2,3,4$,  
 $Dv_5=v_1v_2t_1+t_1^4$, $Dx=0$  and $Dy=x^5+v_1v_3t_1^2$.
They are  contained only in ${\mathcal T}_0(X\times \C P^4)$.

(2) Both ${\mathcal T}_1(X)$ and ${\mathcal T}_1(X\times \C P^4)$
are  same as one   in Theorem 2.2.
Notice that  $P_{0,0}P_{0,1}P_{3,2}P_{2,1}$  is  attached by  a 2-cell
in ${\mathcal T}_2(X\times \C P^4)$
from  $Dv_i=0$ for $i=1,2,3$, $Dv_4=v_1v_2t_1+t_1^4$,
$Dv_5=t_2^4$,  $Dx=0$  and
$Dy=x^5+v_1v_3t_1t_2$.
So ${\mathcal T}(X\times \C P^4)={\mathcal T}(Y)$
for $Y=S^3\times S^3\times S^3\times S^{7}\times S^{9}$.
\hfill\qed

 \begin{rem}
The author must mention about 
the  spaces $X_1$ and $X_2$ in   \cite[Examples 3.8 and 3.9]{Y}
such that ${\mathcal T}_1(X_1)={\mathcal T}_1(X_2)$.
We can check that 2-cells attach on both $P_0P_5P_9P_8$
of them (compare \cite[p.506]{Y}).  
\end{rem}
 
\begin{rem}
In \cite[Question 1.6]{Y},
a rigidity problem is proposed.
It says that  does ${\mathcal T}_0(X)$ with  coordinates   determine ${\mathcal T}_1(X)$ ?.
For  ${\mathcal T}(X)$, it is false as we see in above examples.
But it seems that  there are certain restrictions. 
For example, is ${\mathcal T}_2(X)$ simply connected ?
  \end{rem}

\end{document}